\documentclass[11pt]{amsart}
\usepackage[numbers]{natbib}
\usepackage{amsmath,amsthm,amsfonts,amssymb}
\usepackage[czech,english]{babel}
\usepackage{graphicx}
\usepackage[OT1]{fontenc}
\usepackage{algpseudocode}
\newtheorem{theorem}{Theorem}

\newtheorem{lemma}{Lemma}

\newtheorem{fact}{Fact}
\newtheorem{definition}{Definition}

\newcommand{\bbbn}{\mathbb{N}}

\DeclareGraphicsExtensions{.pdf}

\title{Treeable Graphings Are Local Limits of Finite Graphs}
\thanks{Supported by grant ERCCZ LL-1201 
and CE-ITI, and by the European Associated Laboratory ``Structures in
Combinatorics'' (LEA STRUCO) P202/12/G061}
\author{Lucas Hosseini}
\address{Lucas Hosseini\\
Centre d'Analyse et de Math\'ematiques Sociales (CNRS, UMR 8557)\\
  190-198 avenue de France, 75013 Paris, France}
  \email{lucas.hosseini@gmail.com}
\author{Patrice~Ossona~de~Mendez}
\address{Patrice~Ossona~de~Mendez\\
Centre d'Analyse et de Math\'ematiques Sociales (CNRS, UMR 8557)\\
  190-198 avenue de France, 75013 Paris, France
  and
     Computer Science Institute of Charles University (IUUK)\\
   Malostransk\' e n\' am.25, 11800 Praha 1, Czech Republic}
 \email{pom@ehess.fr}

\date{\today}

\begin{document}
\begin{abstract}
	Let $\mathbf G$ be a graphing, that is a Borel graph defined by $d$ measure preserving involutions. We prove that if $\mathbf G$ is {\em treeable} then it arises as the local limit of some sequence 
	$(G_n)_{n\in\bbbn}$ of graphs with maximum degree at most $d$.
	 This extends a result by Elek \cite{Elek2010} (for $\mathbf G$ a treeing) and consequently extends the domain of the graphings for which Aldous--Lyons conjecture \cite{Aldous2006} is known to be true.
\end{abstract}
\maketitle

\section*{Introduction}

Limits of bounded-degree graph sequences have been extensively studied over the past few years, with the fundamental notion of {\em local convergence} formalized by Benjamini and Schramm \cite{Benjamini2001}. In that framework, a sequence of bounded-degree graphs converges if the distributions of the isomorphism types of its finite radius balls converge, and the limit of a converging sequence is represented by a unimodular probability measure on the space of rooted connected countable graphs or, equivalently, by means of a {\em graphing} (see \cite{Elek2007b} and also  \cite{Aldous2006, Benjamini2001,Gaboriau2005}).
The famous Aldous--Lyons conjecture \cite{Aldous2006} asks whether the converse holds true, that is if every graphing is the local limit of a sequence of finite bounded-degree graphs.

The importance of Aldous--Lyons conjecture goes far beyond the study of graph limits. As mentioned in \cite{Aldous2006}, a positive answer to this question would answer a question of Weiss \cite{Weiss2000}, by showing that all finitely generated groups are “sofic”, and thus would establish several conjectures, since they are known to hold for sofic groups: the direct finiteness conjecture of Kaplansky \cite{kaplansky1972fields} on group algebras, a conjecture of Gottschalk \cite{Gottschalk1973} on “surjunctive” groups in topological dynamics, the Determinant Conjecture on Fuglede-Kadison determinants, and Connes' Embedding Conjecture for group von Neumann algebras \cite{Connes1976}.

Up to now, the cases where the conjecture is known to be true are the hyperfinite case (see e.g. \cite{LovaszBook}) and the case were the graphing is acyclic \cite{Elek2010}, that is when the graphing is a {\em treeing}.

The notions of graphing and treeing being used in different contexts with slightly different definitions, we would like to stress that, in this paper, the graphings and treeings we consider always have uniformly bounded degree, what is usual in the study of local limits of graphs with bounded degrees. 
However, the notion of {\em treeable} Borel equivalence relation --- that is a Borel equivalence relation whose equivalence classes are the connected components of some treeing \cite{kechris2004topics} --- allow graphings and treeing to have countable degree, or locally finite degree \cite{Gaboriau2000}.
 This notion has shown its relevence, as reflected by the proof by Elek and Lippner that every treeable equivalence relation is sofic \cite{Elek20101692}.
 
We say that a graphing $\mathbf G$ is {\em treeable} if there exists a treeing with the same connected components, and that it is {\em weakly treeable} if there exists a treeing $\mathbf T$, such that the domain of $\mathbf G$ is a Borel subset of the domain of $\mathbf T$ with non-null measure, and every connected component of $\mathbf G$ is included in a connected component of $\mathbf T$. 
Obviously, every treeable graphing is weakly treeable.

In this paper we prove (in a purely combinatorial way) that Aldous-Lyons' conjecture holds for weakly treeable graphings.
In order to give an elementary proof of this result, we first give
a proof that colored treeings can be approximated by finite colored graphs (Theorem~\ref{cor:treeing}), and deduce, using model theoretical interpretation techniques, that weakly treeable graphings are limits of finite graphs (Theorem~\ref{thm:treeing-limit}).

\section{Preliminaries}
Let $V$ be a standard Borel space. A graph $\mathbf G$ with vertex set $V$ is a {\em Borel graph} if its edge set is a Borel subset of $V\times V$. Let $d\in\bbbn$ and let $\mathbf G$ be a Borel graph with maximum degree at most $d$. Assume $V$ is endowed an atomless probability measure $\lambda$.
	Then $\mathbf G$ is a {\em graphing} if, for every measurable subsets $A,B$ of $V$ the following {\em Intrinsic Mass Transport Principle} holds:
\begin{equation}
\label{eq:mtp}
\int_A \deg_B(x) \,{\rm d}\lambda(x) = \int_B \deg_A(x) \,{\rm d}\lambda(x),	
\end{equation}
where $\deg_B(x)$ denotes the {\em degree} of $x$ in $B$, that is the number of neighbors of $x$ in $B$.
Above notions extends to vertex/edge colored graphs: for a colored Borel graph it is required that every color class (of vertices or edges) is Borel, while the intrinsic mass transport principle is easily seen to be independent from the coloration.

Several characterizations of graphings are known. For instance, 
a Borel graph $\mathbf G$ with vertex set $V$ and degree at most $d$ is a graphing if and only if  every measurable bijection $T$ with the property that $T(x) = x$ or $T(x)$ is adjacent to $x$ for every $x\in V$ is measure preserving. 
Also, graphings may be characterized as those graphs that can be defined by a finite set $f_1, \dots, f_d$ of measure-preserving involutions on a standard Borel space $V$ by the adjacency of $x$ and $y$ being equivalent to existence of $i\in 1,\dots,d$ such that $f_i(x)=y$. Note that we may assume that for adjacent $x$ and $y$ in $V$ there exists a unique $i\in 1,\dots,d$ such that $f_i(x)=y$. 
 For a more detailed discussion on graphings we refer the reader to \cite{LovaszBook}.

A {\em treeing} is an acyclic graphing.
In other words, a treeing is a graphing, all connected components of which are trees. Treeings have been introduced by Adams \cite{Adams1990} to study Borel equivalence relations.
Recall that in this paper, graphings and treeings are required to have uniformly bounded degrees.
\begin{definition}
 A graphing $\mathbf G$ is {\em treeable}
if there exists a treeing $\mathbf T$ with same vertex set as $\mathbf G$ such that $\mathbf T$ and $\mathbf G$ have the same connected components. 

The graphing $\mathbf G$ is {\em weakly treeable}
if there exists a treeing $\mathbf T$, such that the vertex set as $\mathbf G$ is a Borel subset of the vertex set of $\mathbf T$  with non-null measure, and such that every connected component of $\mathbf G$ is included in a connected component of $\mathbf T$. 
\end{definition}

For a finite graph $G$, we denote by $|G|$ (resp. $\|G\|$) the number of vertices (resp. of edges) of $G$. A {\em rooted} graph is a graph with a distinguished vertex, called the {\em root}. Given a vertex $v$ of a graph $G$ and an integer $r$, the {\em ball} $B_r(G,v)$ of radius $r$ of $G$ centered at $v$ (also called {\em $r$-ball} of $G$ centered at $v$) is the subgraph of $G$ rooted at $v$ induced by all the vertices at distance at most $r$ from $v$. For integer $c$, a {\em $c$-colored} graph
   is a graph with a coloring $\kappa: V\rightarrow\{1,\dots,c\}$ of its vertices in $c$ colors (without any further requirement, like different colors on adjacent vertices); a $c$-colored graphing is a graphing with a measurable coloring $\kappa$
   of its vertices in $c$ colors.

For a integers $d,c,r$ we define $\mathcal G^{d,c}_r$ as the set of all the isomorphism types of $c$-colored rooted graphs with maximum degree at most $d$ and {\em radius} (that is maximum distance of vertices to the root) at most $r$, and define $\mathcal T^{d,c}_r$ as the subset of $\mathcal G^{d,c}_r$ of all the isomorphism types of $c$-colored rooted trees (with maximum degree at most $d$ and radius at most $r$).
	For $k>\ell, t\in \mathcal G^{d,c}_k,$ and $\tau\in \mathcal G^{d,c}_\ell$ we say that $t$ {\em refines} $\tau$ and write $t\prec\tau$ if the $\ell$-ball of the root of (a rooted graph with isomorphism type) $t$ has isomorphism type $\tau$.
	For $k\in\bbbn, t\in\mathcal G_{k+1}^{d,c},$ and $\tau\in\mathcal G_{k}^{d,c}$,  we denote by ${\rm adm}(t,\tau)$ the number of neighbors of the root of (a rooted graph with isomorphism type) $t$ whose $k$-ball has isomorphism type $\tau$. (Note that this does not depend on the choice of a rooted graph with isomorphism type $t$ thus is a well defined value.)

The topological space of finite colored graphs with maximum degree at most $d$ and their local limits (with topology of local convergence) is a compact metrizable space. Its topology is, in particular, induced by the metric
    $\mathrm{dist}(G, H)$ between two colored graphs or graphings $G$ and $H$ 
    defined by
	$$ \mathrm{dist}(G, H)=\sum_{r \geq 1} 2^{-r} \sup_{t\in\mathcal G_{r}^{d,c}} | \Pr(B_r(G,X)\simeq t) - \Pr(B_r(H,Y)\simeq t)|, $$
	where $X$ (resp. $Y$) stands for a uniformly distributed random vertex of $G$ (resp. of $H$).
	
	It is easily seen that this defines an ultrametric on graphings.

	In this paper, we shall make use of the following result of 
	Kechris, Solecki, and Todorcevic:
\begin{lemma}[\cite{Kechris1999}]
	\label{lem:kechris}
	Let  $G$ be a Borel graph with maximum degree at most $d$. Then $G$ has a proper vertex Borel coloring with at most $d+1$ colors and a proper edge Borel coloring with at most $2d-1$ colors.
\end{lemma}

An easy generalization of this lemma to rainbow coloring of bounded radius balls will be particularly useful, when used together with the model theoretical notion of interpretation. In order to prove this extension, we first prove that the $k$-th power $\mathbf G^k$ of a graphing $\mathbf G$ is a graphing.
\begin{lemma}
	\label{lem:power}
	Let $\mathbf G$ be a graphing and let $r\in\bbbn$. Then the $r$-th power $\mathbf G^r$ of $\mathbf G$ is a graphing.
\end{lemma}
\begin{proof}
	Let $Z$ be the set of all sequences $(i_1,\dots,i_k)$ of integers in $1,\dots,d$ with length $k\leq r$, ordered by lexicographic order, together with a special element $\omega$. Let $V$ be the vertex set of $G$. For $x\in V$ define $z_x:V\rightarrow Z$ by
	$$z_x(y)=\begin{cases}
	\omega&\text{ if }{\rm dist}_{\mathbf G}(x,y)>r\\
   \min\{(i_1,\dots,i_k)\in Z:\ y=f_{i_1}\circ\dots\circ f_{i_k}(x)\}&\text{otherwise}	
\end{cases}$$
	Let $A,B$ be measurable subsets of $V$. Then in $\mathbf G^r$ it holds
\begin{align*}
\int_A&{\rm deg}_B(x)\,{\rm d}\lambda(x)\\
&=\sum_{(i_1,\dots,i_k)\in Z\setminus\{\omega\}} \lambda(\{x\in A: \exists y\in B\ z_x(y)=(i_1,\dots,i_k))\\
&=		\sum_{(i_1,\dots,i_k)\in Z\setminus\{\omega\}} \lambda(f_{i_1}\circ\dots\circ f_{i_k}(\{y\in B:  \exists x\in A\ z_y(x)=(i_k,\dots,i_1)\}))\\
	&=
		\sum_{(i_1,\dots,i_k)\in Z\setminus\{\omega\}} \lambda(\{y\in B: \exists x\in A\ z_y(x)=(i_k,\dots,i_1)\})\\
	&=\int_B{\rm deg}_A(y)\,{\rm d}\lambda(y)
\end{align*}
Thus the intrinsic mass transport holds for $\mathbf G^r$.
\end{proof}

\begin{lemma}
	\label{lem:coloring-balls}
	Let $\mathbf G$ be a graphing, and let $r \in \mathbb{N}$. Then there exists a Borel coloring of $\mathbf G$ with finitely many colors such that for every $x \in G$, the induced coloring on $B_r(\mathbf G, x)$ is rainbow (i.e. no two vertices at distance at most $2r$ have the same color).
\end{lemma}
\begin{proof}
	Let $\mathbf G^{2r}$ be the graphing on the same probability space as $G$ where $x$ and $y$ are adjacent in $\mathbf G^{2r}$ if ${\rm dist}_{\mathbf G}(x, y) \leq 2r$ (this is a graphing according to Lemma~\ref{lem:power}). 
	According to Lemma~\ref{lem:kechris}, $\mathbf G^{2r}$ has a proper Borel vertex coloring with finitely many colors. This coloring, considered as a coloring of $\mathbf G$, satisfies the requirements of the lemma.
\end{proof}

Recall that a first-order formula $\phi$ (in the language of colored graphs) is {\em local} if its satisfaction only depends on a fixed neighborhood of its free variables, and it is {\em strongly-local} if $\phi$ is local and there exists $r \in \mathbb{N}$ s.t. the satisfaction of $\phi$ implies that all its free variables are pairwise at distance at most $r$.
For instance, the formula $\phi_1(x_1,x_2)$ asserting that $x_1$ and $x_2$ are at distance at most $100$ is local, although the formula $\phi_2(x_1)$ asserting that the graph contains a vertex of degree $1$ not adjacent to $x_1$ is not.

For a graph $G$ and a formula $\phi$ with $p$ free variables, we denote by $\phi(G)$ the set of satisfying tuples of $\phi$ for $G$:
$$\phi(G)=\{(v_1,\dots,v_p)\in G^p:\ G\models \phi(v_1,\dots,v_p)\},$$
and we denote by $\langle\phi,G\rangle$ the probability that $\phi$ is satisfied for a random assignment of vertices of $G$ to the free variables of $\phi$:
$\langle\phi,G\rangle:=|\phi(G)|/|G|^p$.

The interest of considering local formulas in the context of local convergence appears clearly from the next lemmas.

\begin{lemma}[\cite{CMUC}]
	\label{lem:loclim}
	Let $(G_n)_{n\in\bbbn}$ be a sequence of finite (colored) graphs with maximum degree at most $d$ and such that $\lim_{n\rightarrow\infty}|G_n|=\infty$. Then the following conditions are equivalent:
	\begin{enumerate}
		\item The sequence $(G_n)_{n\in\bbbn}$ is local convergent;
		\item for every local formula $\phi$ with a single free variable, the probability
$\langle\phi,G_n\rangle$ that $\phi$ is satisfied for a random assignment of the free variable converges as $n\rightarrow\infty$;
		\item for every local formula $\phi$, the probability
$\langle\phi,G_n\rangle$ that $\phi$ is satisfied for a random assignment of the free variables converges as $n\rightarrow\infty$.	
\end{enumerate}
\end{lemma}

In this paper, an {\em interpretation scheme} $\mathsf{I}$ will be
defined by some first-order formulas $\xi,\eta,\mu_1,\dots,\mu_q$, where $\eta$ has two free variables, and is symmetric and anti-reflexive, and where $\xi$ and the $\mu_i$'s have a single free variable. Moreover it is required that $\eta(x,y)\rightarrow\xi(x)\wedge\xi(y)$ and $\mu_i(x)\rightarrow\xi(x)$.

An interpretation scheme  $\mathsf I$ allows to define a new colored (or marked) graph $\mathsf I(G)$ from a (colored) graph $G$ as follows: the vertex set of $\mathsf I(G)$ is the set $\xi(G)$ of vertices of $G$ that satisfy $\xi$; the edge set of $\mathsf I(G)$ is the set $\eta(G)$ of pairs $\{x,y\}$ of vertices of $G$ that satisfy $\eta(x,y)$, and the set of vertices of $\mathsf I(G)$ marked with mark $i$ is $\mu_i(G)$.
An interpretation scheme $\mathsf I$ is {\em local} (resp. {\em strongly local}) if all the formulas $\xi,\eta,\mu_1,\dots,\mu_p$ are local (resp. strongly local). A basic fact about
local interpretation scheme is the following \cite{Hodges1997}:
\begin{fact}
\label{fact:1}
	Let $\mathsf I$ be a local interpretation scheme. Then for every local formula $\phi$ there is a local formula $\mathsf I^\ast(\phi)$ such that for every colored graph $G$ it holds
	$$
	{\mathsf I}^\ast(\phi)(G)=\phi(\mathsf I(G)).
	$$
\end{fact}

A particular case of local interpretation is taking the $k$-th power. Hence the following lemma is a generalization of Lemma~\ref{lem:power}.
\begin{lemma}
\label{lem:graphing_interpretation}
	Let $\mathbf G$ be a graphing (with associate probability measure $\nu$) and let $\mathsf I=(\xi,\eta,\mu_1,\dots,\mu_p)$ be a strongly local interpretation scheme.
	
	Assume $\nu(\xi(\mathbf G))>0$ and let $\widehat\nu$ be the probability measure on $\xi(\mathbf G)$ defined by $\widehat\nu(X)=\nu(X)/\nu(\xi(\mathbf G))$.
	
	Then $\mathsf I(\mathbf G)$, with associate probability measure
	$\widehat\nu$, is a graphing.
\end{lemma}
\begin{proof}
Let $V$ be the vertex set of $\mathbf G$.
Then there exists an integer $r$ such that the satisfaction of any of the formulas $\xi,\eta,\mu_1,\dots,\mu_p$ only depends on the $r$-neighborhood of the first free variable and implies that the other free variables (if any) are at distance at most $r$ from the first one.
According to Lemma~\ref{lem:coloring-balls} we can assume that the graphing $\mathbf G$ is colored in such a way that no two vertices with same color exist in a same ball of radius $r$. 
 For each $t\in\mathcal G_{r}^{d,c}$ we fix a rooted graph $(F_t,r_t)$ with isomorphism type $t$. Thanks to the strong locality of $\eta$, there exists a finite family $\mathcal F$ of pairs $(t,z)$ with $t\in\mathcal G_{2r}^{d,c}$ and $z\in V(F_t)$, such that for every two vertices $u,v\in V$ it holds that $\mathbf G$ satisfies $\eta(u,v)$ if and only if
there is some $(t,z)\in\mathcal{F}$ such that there is a (necessarily unique) isomorphism  $B_{r}(\mathbf G,x)\rightarrow F_t$ mapping $x$ to $r_t$ and $y$ to $z$. 
Let $g_{t,z}$ be the map defined as follows: Let $x\in V$.
\begin{itemize}
\item if $B_{r}(\mathbf G,x)\simeq t$ then $g_{t,z}(x)$ is the unique vertex $y$ in $B_{2r}(\mathbf G,x)$ with same color as $z$;
\item if $x$ has the same color as $z$ and if there exists a vertex $y$ at distance at most $2r$ from $x$ such that
$B_{r}(\mathbf G,y)\simeq t$ then $g_{t,z}(y)=x$;
\item otherwise, $g_{t,z}(x)=x$.
\end{itemize}
 According to the intrinsic mass transport principle, $g_{t,z}$ is a measure preserving involution, and it is immediate that the edges of $\mathsf{I}(\mathbf{G})$ are those pairs $\{x,y\}$ such that $x\neq y$ and there exists $(t,z)\in\mathcal F$ with $y=g_{z,t}(x)$. Also, as the satisfaction of $\mu_i(x)$  only depends on the isomorphism type of $B_r(\mathbf G,x)$ the set $\mu_i(\mathsf{I}(\mathbf{G}))$ is Borel. It follows that $\mathsf{I}(\mathbf{G})$ is a graphing.
\end{proof}

From Fact~\ref{fact:1} and Lemmas~\ref{lem:loclim} and~\ref{lem:graphing_interpretation} we immediately deduce:
\begin{lemma}
	\label{lem:interpretation-continuous}
	Let $(G_n)_{n\in\bbbn}$ be a local convergent sequence of finite (colored) graphs with maximum degree at most $d$ such that $\lim_{n\rightarrow\infty}|G_n|=\infty$, and let $\mathsf I=(\xi,\eta,\mu_1,\dots,\mu_p)$ be a strongly local interpretation scheme with 
	$$\lim_{n\rightarrow\infty}\xi(G_n)>0.$$
	
	Then the sequence $(\mathsf I(G_n))_{n\in\bbbn}$ is local convergent. Moreover, if the sequence $(G_n)$ converges locally towards a graphing $\mathbf G$, then $\mathsf I(G_n)$ converges locally towards the graphing $\mathsf I(\mathbf{G})$.
	\qed
\end{lemma}

In other words, a local interpretation scheme defines a continuous map.

\section{Finite Approximations of Treeings}

Let $k\in\bbbn$ and let $\epsilon>0$ be a positive real.
We now describe how a treeing $\mathbf T$ can be approximated by a finite graph.

Let $\mathbf{T}$ be a $c$-colored treeing (with associated probability distribution $q$ on $\mathcal{G}_{k+1}^{d,c}$). Without loss of generality, we can assume according to Lemma~\ref{lem:coloring-balls} that $c$ induces a proper coloring on $B_k(\mathbf T, x)$ for all $x \in \mathbf T$. It follows that any type $t \in \mathcal{T}_{k+1}^{d,c}$ that is realized as a $k$-ball around a vertex of $\mathbf T$ is such that for any $\tau \in \mathcal{T}_{k}^{d,c}$, $\mathrm{adm}(t, \tau) \in \{ 0, 1 \}$.

Let $N(k,d,c,\epsilon) = \frac{d^{3k+3}|\mathcal T_k^{d,c}|^2+d^{k+1} |\mathcal T_k^{d,c}|\,|\mathcal T_{k+1}^{d,c}|}{\epsilon}$ and let $n > N(k,d,c,\epsilon)$.

We define a finite graph $G$ as follows.
	
For $t \in \mathcal{T}_{k+1}^{d,c}$ let $V_t$ be a set of cardinality $\lceil n q(t) \rceil$, and let $V$ be the disjoint union of the $V_t$. 
For $\tau\in \mathcal T_k^{d,c}$ it will be convenient to denote by  
$W_\tau$ the union of the $V_t$ such that $t\prec \tau$.
	
Let $E$ be a maximal subset of pairs of elements of $V$ such that the following properties hold in the graph $G = (V, E)$:
\begin{enumerate}
	\item the girth of $G$ is at least $2k+4$,
	\item no vertex in $V_t$ has more than $\mathrm{adm}(t,\tau)$ neighbors in $W_\tau$ (for every $t \in \mathcal{T}_{k+1}^{d,c}$ and $\tau \in \mathcal{T}_k^{d,c}$).
\end{enumerate}

Accordingly, we say that $x \in V_t$ is {\em $\tau$-deficient} if $\deg_{W_{\tau}}(x) < \mathrm{adm}(t, \tau)$. 
A vertex is {\em bad} if it is $\tau$-deficient for some $\tau \in \mathcal{T}_k^{d,c}$, and it is {\em perfect} if no vertex at distance at most $k$ from it is bad.

In order to prove that $G$ is an approximation of $\mathbf{T}$ we will need the following two lemmas.

\begin{lemma}
\label{lem:nonperfect-bound}
The number of non-perfect vertices of $G$ is at most $d^{3k+3}|\mathcal T_k^{d,c}|^2+d^{k+1} |\mathcal T_k^{d,c}|\,|\mathcal T_{k+1}^{d,c}|$.
\end{lemma}
\begin{proof}
	We first bound the number of $\tau'$-deficient vertices in $W_\tau$ for $\tau,\tau'\in\mathcal T_k^{d,c}$. (Note that since we consider graphs with girth at least $2k+4$, the $k+1$-ball of every vertex is a rooted tree.)

If there exists $x \in W_{\tau'}$ that is $\tau$-deficient then any $y \in W_{\tau}$ that is $\tau'$-deficient is at distance at most $2k+3$ from $x$ (by maximality of $E$). Hence there are at most $d^{2k+3}$ vertices in $W_{\tau}$ that are $\tau'$-deficient.

Otherwise, for every $x\in V_{t'}\subseteq W_\tau'$ it holds $\deg_{W_{\tau}}(x) = \mathrm{adm}(t', \tau)$. Hence
$$\sum_{t \prec \tau} \sum_{y \in V_t} \deg_{W_{\tau'}}(x)=\sum_{t' \prec \tau'} \sum_{x \in V_{t'}} \deg_{W_{\tau}}(x)=\sum_{t' \prec \tau'} |V_{t'}|\, \mathrm{adm}(t', \tau)$$
(first equality follows from double-count of the edges between $W_\tau$ and $W_{\tau'}$). For $t \in \mathcal{T}^{d,c}_{k+1}$, define $q(t)=\Pr[B_r(T,X) \simeq t]$. Applying the mass transport principle to
 $A = \{ x ; B_k(\mathbf{T}, x) \simeq \tau\}$ and $B = \{ y ; B_k(\mathbf{T}, y) \simeq \tau' \}$ we get
		$$\sum_{t \prec \tau} q(t)\, \mathrm{adm}(t, \tau')=\sum_{t' \prec \tau'} q(t')\, \mathrm{adm}(t', \tau).$$
		Hence  the number of $\tau'$-deficient vertices in $W_\tau$ is bounded by
\begin{align*}
\sum_{t \prec \tau} \sum_{x \in V_{t}}&(\mathrm{adm}(t, \tau')-\deg_{W_{\tau'}}(x))\\
	&\leq \sum_{t \prec \tau} |V_{t}|\, \mathrm{adm}(t, \tau')-\sum_{t' \prec \tau'} |V_{t'}|\, \mathrm{adm}(t', \tau)\\
	&\leq \sum_{t \prec \tau} (Nq(t)+1)\, \mathrm{adm}(t, \tau')-\sum_{t' \prec \tau'} Nq(t')\, \mathrm{adm}(t', \tau)\\
	&\leq dC(\tau),
\end{align*}
where $C(\tau)=|\{t\in\mathcal T_{k+1}^{d,c}: t\prec\tau\}|$. 

Altogether, the number of of $\tau'$-deficient vertices in $W_\tau$ is bounded by $\max(d^{2k+3},dC(\tau))<d^{2k}+dC(\tau)$.
Summing up over pairs $\tau,\tau'$, we get that the number of bad vertices of $G$ is bounded by 
	$$\sum_{\tau}\sum_{\tau'}(d^{2k+3}+dC(\tau))=d^{2k+3}|\mathcal T_k^{d,c}|^2+d|\mathcal T_k^{d,c}|\,|\mathcal T_{k+1}^{d,c}|.$$
		Since a non-perfect vertex is one that is at distance less or equal than $k$ from a bad vertex, the number of non-perfect vertices is bounded by $d^{3k+3}|\mathcal T_k^{d,c}|^2+d^{k+1} |\mathcal T_k^{d,c}|\,|\mathcal T_{k+1}^{d,c}|$.
\end{proof}

\begin{lemma}
\label{lem:perfect-type}
If $x \in W_\tau$ is a perfect vertex of $G$ then $x$ has type $\tau$.
\end{lemma}
\begin{proof}
We prove by induction that for every $t \in \mathcal T_{k+1}^{d,c}$ and every vertex $x \in V_t$, if every vertex at distance at most $r$ from $x$ is $good$, then $B_r(G, x) \simeq t_r$ where $t_r$ is the $r$-ball of $t$ centered on its root.

This trivially holds for $r = 0$ by construction of $G$.

Suppose it holds for $r$. Let $x$ be a vertex of $(r+1)$-type $t$ s.t. every vertex at distance at most $r+1$ from $x$ is good.

Let $y_1, \dots, y_p$ be the neighbords of $x$ in $G$, $s_1, \dots, s_p$ the neighbors of the root $r$ of $t$ (of respective $r$-types $\tau_1, \dots, \tau_p$), and let $\tau$ be an $r$-type.

By construction $|\{ i ; \tau_i \simeq \tau \}| = \mathrm{adm}(t, \tau)$. Moreover, by induction, $|\{ i ; \tau_i \simeq \tau \}| = |\{ i ; B_r(G, y_i) \simeq \tau \}|$.

Hence, we can assume $B_r(G, y_i) \simeq \tau_i$.

Thus, for each $i$ there exists an isomorphism $f_i : B_r(G, y_i) \to \tau_i$ s.t. $y_i \mapsto s_i$. 
By color-preservation, $f_i(x) = r$.

For each $i$, let $T_i$ (resp. $T'_i$) be the subgraph of $G$ (resp. $t$) induced by the vertices of the 
connected component of $y_i$ (resp. $s_i$) in $B_r(G, x)$ (resp. $t$) when removing the edge 
$\{x, y_i\}$ (resp. $\{r, s_i\}$).

For each $i$, $f_i|_{T_i}$ is an isomorphism between $(T_i, y_i)$ and $(T'_i, s_i)$.

Finally, we define $f : B_{r+1}(x) \to t$ by $f(x) = r$ and $f(y) = f_i(y)$ if $y \in T_i$, which is indeed an isomorphism between $B_{r+1}(x)$ and $t$.
\end{proof}

\begin{lemma}
\label{lem:treeing}
	Let $\mathbf{T}$ be a $c$-colored treeing with maximum degree $d$.
	
	Then, for every integer $r$ and every positive real $\epsilon>0$ there exists integer $n_0$ such that for every $n\geq n_0$ there exists a $c$-colored finite graph $G$ with maximum degree $d$ and satisfying the following properties:
	\begin{enumerate}
		\item  $n\leq |G| \leq n+|\mathcal T_{r+1}^{d,c}|$, 
		\item $\mathrm{girth}(G) \geq 2r+1$,
		\item for every $t \in \mathcal{G}_r^{d,c}$ it holds
 $$\bigl|\Pr[B_r(G, X) \simeq t] - \Pr[B_r(\mathbf T, X) \simeq t]\bigr| < \epsilon.$$
	\end{enumerate}
	 In particular,
 ${\rm dist}(G,\mathbf T)<\epsilon+2^{-r}$.
\end{lemma}

\begin{proof}
	By construction $\mathrm{girth}(G) \geq 2r+1$ and $n \leq |G| \leq n+|\mathcal{T}_{r+1}^{d,c}|$.
	
	Moreover, following Lemmas~\ref{lem:nonperfect-bound} and \ref{lem:perfect-type} the number of vertices $x \in W_\tau$ such that $B_r(G, x) \not\simeq \tau$ is bounded by $d^{r}$ times the number of non-perfect vertices. Hence it holds for every $\tau \in \mathcal{G}_r^{d,c}$:
	$$|\Pr( B_r(G, X) \simeq \tau) - \Pr( B_r(\mathbf T, X) \simeq \tau)| \leq \frac{d^{3r+3}|\mathcal T_r^{d,c}|^2+d^{r+1} |\mathcal T_r^{d,c}|\,|\mathcal T_{r+1}^{d,c}|}{n} < \epsilon.$$
\end{proof}

As a corollary we have an alternative proof of the following result \cite{Elek2010}:
\begin{theorem}
	\label{cor:treeing}
	Every $c$-colored treeing with maximum degree $d$ is the local limit of a sequence of finite $c$-colored graphs with maximum degree $d$ and girth growing to infinity.
\end{theorem}
\begin{proof}
	For $n\in\bbbn$. According to Lemma~\ref{lem:treeing} there exists a $c$-colored graph $G_n$ with maximum degree at most $d$, order at least $n$, girth at least $2n+1$, and such that 
	$$\sup_{t \in \mathcal{G}_n^{d,c}}
 \bigl|\Pr[B_n(G_n, X) \simeq t] - \Pr[B_n(\mathbf T, X) \simeq t]\bigr| < 2^{-n}/|\mathcal{G}_n^{d,c}|.$$
 Hence for every integer $r$ and every $n\geq r$ it holds
 $$\sum_{t \in \mathcal{G}_r^{d,c}}\bigl|\Pr[B_r(G_n, X) \simeq t] - \Pr[B_r(\mathbf T, X) \simeq t]\bigr| < 2^{-n}.$$
 
 It follows that $\mathbf T$ is the local limit of the sequence
 $(G_n)_{n\in\bbbn}$.
\end{proof}

\section{Finite Approximations of Weakly Treeable Graphings}
Finite approximation of weakly treeable graphings will makes use of the model theoretical notion of interpretation.

\begin{lemma}
	\label{lem:neighborhood-inclusion}
	Let $\mathbf G$ be a weakly treeable graphing and $\mathbf T$ be a treeing, such that the domain of $\mathbf G$ is a Borel subset of the domain of $\mathbf T$  with non-null measure, and each connected component of $G$ is included i a connected component of $\mathbf T$.
	Then it holds:
		$$\forall \epsilon > 0, \exists r, \nu(\{x; B_1(\mathbf G, x) \subset B_r(\mathbf T, x)) > 1 - \epsilon.$$
\end{lemma}
\begin{proof}
	Let $X_r = \{ x; B_1(\mathbf{G}, x) \subseteq B_r(\mathbf{T}, x)\}$. Each $X_r$ is measurable as it is definable in first order logic, and the sequence $(X_r)$ is monotonous. Hence, according to Beppo-Levi, $\nu(\cup_{r} X_r) = \lim_{r \to \infty} \nu(X_r)$, and since $\cup_r X_r = G$, $\nu(\cup_r X_r) = 1$, so for all $\epsilon$, there is an $r$ s.t. $\nu(X_r) > 1 - \epsilon$.
\end{proof}

\begin{lemma}
	\label{lem:measure-bound}
	Let $\mathbf G, \mathbf H$ be graphings of degree bounded by $d$ on the same domain s.t. $\nu(\{ x ; B_1(\mathbf{G}, x) \neq B_1(\mathbf{H}, x) \}) \leq \epsilon$ and $r \in \mathbb{N}$. 
	Then $\nu(\{ x ; B_r(\mathbf{G}, x) \not\simeq B_r(\mathbf{H}, x) \}) < d^r \epsilon$.
\end{lemma}
\begin{proof}
	Let $X \subset G$. Then $\nu(B_1(\mathbf{G}, X)) \leq d \nu(X)$ by using a caracterization of $G$ in terms of $d$ measure-preserving involutions. 
	
	Moreover, 
	if $B_r(\mathbf{G}, x) \not\simeq B_r(\mathbf{H}, x)$ then there exists $y \in B_r(\mathbf{G}, x)$ s.t. $B_1(\mathbf{G}, y) \neq B_1(\mathbf{H}, y)$.
	 
	Hence, $\nu(\{ x ; B_r(\mathbf{G}, x) \neq B_r(\mathbf{H}, x)\}) < d^r \epsilon$.
\end{proof}

\begin{lemma}
	\label{lem:treeing-sequence}
	Let $\mathbf G$ be a weakly treeable graphing. There exists a sequence $(H_n)_n$ of finite graphs with a countable marking s.t.
	\begin{itemize}
		\item $\mathrm{girth}(H_n) \to \infty,$
		\item $\forall \epsilon > 0, \exists r, \limsup_{n \to \infty} \mathrm{dist}(\mathbf G, I_r(H_n)) < \epsilon.$
	\end{itemize}
\end{lemma}
\begin{proof}
	Let $\mathbf T$ be a treeing such that the domain of $\mathbf G$ is a Borel subset with positive measure of the domain of $\mathbf T$, and such that every connected component of $\mathbf G$ is included in a connected component of $\mathbf T$. We extend the coloring of $\mathbf T$ by adding a mark $\xi$ such that $\xi(x)$ holds exactly if $x\in\mathbf G$.
	
	Let $\epsilon > 0$. According to Lemma~\ref{lem:neighborhood-inclusion}, there is $r$ s.t. $\nu(\{x; B_1(\mathbf G, x) \subset B_r(\mathbf T, x)\}) > 1 - \epsilon$. Let $c$ be an $r$-ball coloring of $\mathbf G$.
	
	From $c$, we define a new coloring $c'$ the following way: for each vertex $x$ of $\mathbf G$, let $u_1, \dots, u_k$ be the vertices of $B_1(\mathbf G, x) \cap B_r(\mathbf{T}, x)$ (i.e. the neighbors of $x$ in $\mathbf G$ that are at distance at most $r$ in $\mathbf T$), we define $c'(x)$ as $(c(x), n(x))$ where $n(x) = \{ c(u_1), \dots, c(u_k)\}$.
	
	We define the interpretation $I_r=(\xi,\eta_r,\dots)$, which preserves the coloring $c$, and such that 
	$\eta_r(x,y)$ holds exactly if $y \in B_r(H, x)$ and $\chi(y) \in n(x)$.
	
	It holds $\mathrm{dist}(\mathbf G, I_r(T)) < d^r\epsilon$ following Lemma~\ref{lem:measure-bound} since the neighborhoods of each vertex are the same in $\mathbf G$ and in $I_r(\mathbf T)$ except for a set of vertices of measure less than $\epsilon$.
	
	Let $(H_n)$ be a sequence of colored finite graphs converging towards $(\mathbf T, c')$. Then $(I_r(H_n))_n$ converges locally towards $I_r(\mathbf T)$ according to Lemma~\ref{lem:interpretation-continuous}.
\end{proof}

\begin{theorem}
	\label{thm:treeing-limit}
	Let $\mathbf G$ be a vertex-colored weakly treeable graphing. Then there exists a sequence $(G_n)$ of finite vertex-colored graphs that converges towards $\mathbf G$. 
\end{theorem}
\begin{proof}
	For $k \in \mathbb{N}$, there exists $r_k, n_k$ s.t. $\mathrm{dist}(\mathbf G, I_{r_k}(H_{n_k})) < 2^{-k}$ according to Lemma~\ref{lem:treeing-sequence}, so we can let $G_k = I_{r_k}(H_{n_k})$.
\end{proof}

\providecommand{\noopsort}[1]{}\providecommand{\noopsort}[1]{}
\providecommand{\bysame}{\leavevmode\hbox to3em{\hrulefill}\thinspace}
\providecommand{\MR}{\relax\ifhmode\unskip\space\fi MR }
\providecommand{\MRhref}[2]{%
  \href{http://www.ams.org/mathscinet-getitem?mr=#1}{#2}
}
\providecommand{\href}[2]{#2}

\end{document}